\documentclass[preprint,10pt,authoryear]{elsarticle}
\newif\ifanon
\ifdefined\ANON\anontrue\else\anonfalse\fi
\usepackage[T1]{fontenc}
\usepackage[utf8]{inputenc}
\usepackage{mathptmx}
\usepackage{microtype}
\usepackage[margin=2.5cm]{geometry}
\usepackage{amsmath,amssymb}
\usepackage{graphicx}
\usepackage{adjustbox}
\usepackage{booktabs}
\usepackage{multirow}
\usepackage{xcolor}
\usepackage{colortbl}
\usepackage{siunitx}
\usepackage{enumitem}
\usepackage{caption}
\usepackage{subcaption}
\usepackage{float}
\usepackage[section]{placeins}
\usepackage{setspace}
\usepackage[colorlinks=true,linkcolor=zblue,citecolor=zblue,urlcolor=zblue]{hyperref}
\usepackage[capitalise,noabbrev]{cleveref}
\usepackage{orcidlink}

\definecolor{zblue}{HTML}{1C5CAB}
\definecolor{zbluebg}{HTML}{EAF2FC}
\definecolor{zgrey}{HTML}{6E6D68}
\setlist{itemsep=1pt,topsep=3pt,leftmargin=*}
\newcommand{\zfmn}{\textsc{Z-FMN}}
\newcommand{\bench}{\textsc{RAGTruth-FMN}}
\newcommand{\ci}[2]{\,{\color{zgrey}\scriptsize[#1,\,#2]}}
\newcommand{\best}[1]{\textbf{#1}}
\newcommand{\ours}{\rowcolor{zbluebg}}
\makeatletter
\def\ps@pprintTitle{\let\@oddhead\@empty\let\@evenhead\@empty\def\@oddfoot{\hfil\thepage\hfil}\let\@evenfoot\@oddfoot}
\makeatother

\begin{document}

\begingroup\let\hrule\relax
\begin{frontmatter}

\title{Hallucination detection in retrieval-augmented language models by Z-number evidence fusion on concept networks}

\ifanon
\author{Anonymous author(s)}
\else
\author[iitj]{Sweta Kaman\corref{cor1}\,\orcidlink{0000-0002-3263-8282}}
\ead{kaman.1@iitj.ac.in}
\cortext[cor1]{Corresponding author. ORCID: \href{https://orcid.org/0000-0002-3263-8282}{0000-0002-3263-8282}}
\affiliation[iitj]{organization={School of AI and Data Science, Indian Institute of Technology Jodhpur},
            city={Jodhpur},
            postcode={343020},
            country={India}}
\fi

\begin{abstract}
Retrieval-augmented large language models are increasingly used in engineering and industrial information systems, yet they still produce hallucinations: statements that their retrieved sources do not support. Detecting these faults automatically is a verification and validation problem for artificial intelligence software. \emph{Contribution in artificial intelligence:} we propose the Z-number forma mentis network (\zfmn{}), a training-free approximate-reasoning detector that turns each answer into a network of concepts and syntactic links and attaches to every concept and link a Z-number whose reliability part is a subjective-logic opinion. Cumulative fusion combines three evidence sources, the generator's confidence across resampled answers, natural language inference against the source and consistency with the resamples, so that each source counts in proportion to the evidence it carries. \emph{Application in engineering:} \zfmn{} acts as a claim-level fault monitor that ranks the links of an answer for human review. On all 1,350 test answers of three open models in the public RAGTruth dataset, it detects hallucinated links with an area under the receiver operating characteristic curve of 0.793, ahead of SelfCheckGPT by 0.068 and of source entailment by 0.026, and within 0.009 of a supervised model, without any training. A reviewer who checks the riskiest 20\% of links finds 60\% of the hallucinated ones, against 50\% with SelfCheckGPT. The uncertainty of the fused opinion also separates contradicted from invented content. We release a concept- and link-level benchmark of 273,456 labelled elements with all scores and code.
\end{abstract}


\end{frontmatter}
\endgroup

\section{Introduction}
\label{sec:intro}

Retrieval-augmented generation (RAG; \citealp{lewis2020retrieval}) grounds a large language model (LLM) in retrieved passages. It has become a standard architecture for engineering and industrial information systems: assistants that answer questions over technical documentation, tools that summarise reports, and pipelines that turn structured records into text. Yet the model still adds facts the passages do not contain and contradicts facts they do. RAGTruth \citep{niu-etal-2024-ragtruth} annotated such spans in about 18,000 answers and found them in a large share of answers from every model it studied.

For an engineer who deploys such a system, a hallucination is a fault in the output of an AI software component. Like any fault, it has to be detected, located and handed to a person or a downstream check before the output is used. Detection alone is not enough: a monitor that flags a whole answer or sentence tells the reviewer that something is wrong, but not which claim to check, and review time is the scarce resource. This paper treats hallucination detection as claim-level fault monitoring and asks how much of the fault a reviewer finds for a given review effort.

Existing detectors fall into three families. \emph{Internal} detectors read the generator's token probabilities, entropy or hidden states \citep{malinin2021uncertainty,kadavath2022language,azaria-mitchell-2023-internal}. \emph{Sampling} detectors such as SelfCheckGPT \citep{manakul-etal-2023-selfcheckgpt} and semantic entropy \citep{kuhn2023semantic,farquhar2024detecting} compare an answer with resampled answers. \emph{Grounding} detectors check entailment against the source \citep{laban-etal-2022-summac,honovich-etal-2022-true,zha-etal-2023-alignscore}. All three return one number per token or per sentence.

This framing misses two things. First, the unit of error is usually a relation. ``The fire was extinguished by the crew'' can be wrong because the wrong agent is attached to the right event, while every word in it appears in the source. Token scores spread such an error over many tokens, and sentence scores blur it with the correct claims around it. Second, the evidence sources disagree \emph{in kind}. An entailment model that returns ``neutral'' has found no evidence either way; one that returns ``contradiction'' has found evidence against. A resample that omits a claim is weak evidence; one that states the opposite is strong. Averaging these signals into a probability discards how much evidence each one actually carries.

\begin{figure}[!htbp]
  \centering
  \includegraphics[width=\textwidth]{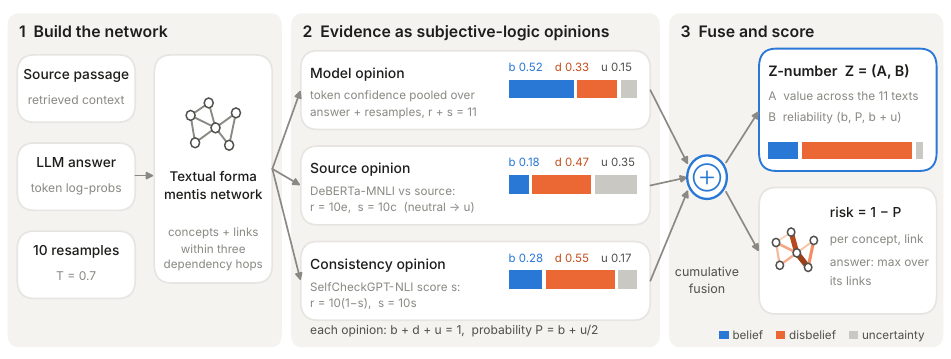}
  \caption{\textbf{\zfmn{} at a glance.} (1) The answer becomes a network of lemmatised concepts linked within three dependency hops. (2) Each concept and link receives three subjective-logic opinions, shown as belief\,/\,disbelief\,/\,uncertainty bars: the generator's confidence pooled over the answer and its resamples, entailment against the source (a neutral verdict adds uncertainty, not disbelief), and SelfCheckGPT-style consistency. (3) Cumulative fusion weights each opinion by its evidence; the fused opinion is the reliability $B$ of a Z-number, and $1-P$ is the element's hallucination risk. Bar values are illustrative.}
  \label{fig:overview}
\end{figure}

We address both with \zfmn{} (\cref{fig:overview}). It builds a textual forma mentis network (TFMN; \citealp{stella2020text}) of each answer, whose nodes are concepts and whose edges are syntactic links between them. Every node and edge receives a hybrid Z-number \citep{zadeh2011note}, a pair of a value and its reliability, whose reliability is a subjective-logic opinion \citep{josang2016subjective} of belief, disbelief and uncertainty. Three evidence sources become three opinions, and cumulative fusion combines them so that a source with much evidence dominates one with little. Nothing is trained, so the monitor can be attached to a new model or domain without labelled data.

\paragraph{Contributions} The contribution in artificial intelligence is an approximate-reasoning method for fusing heterogeneous evidence on a graph; the engineering application is claim-level fault monitoring of retrieval-augmented systems.
\begin{enumerate}
  \item \textbf{\zfmn{}}, a training-free hallucination detector that scores the concepts and links of an answer by fusing model, source and consistency evidence as Z-number opinions (\cref{sec:method}).
  \item \textbf{A deployment analysis} of \zfmn{} as a fault monitor: how many hallucinated claims a reviewer finds when checking only the riskiest links, what each answer costs, and how to set the operating point (\cref{sec:deploy}).
  \item \textbf{\bench{}}, a concept- and link-level benchmark with 273,456 labelled elements from 1,350 answers of three open models, derived from RAGTruth's human span annotations and released with every score in this paper (\cref{sec:bench}).
  \item \textbf{A state-of-the-art training-free result.} \zfmn{} beats SelfCheckGPT-NLI on links by 0.068 and on concepts by 0.033, wins against it on all nine model-and-task datasets, and comes within 0.009 of a supervised model (\cref{sec:main}).
  \item \textbf{An analysis of what the Z-number adds.} Its reliability component gives a significant gain, its uncertainty dimension tells contradictions from invented content, and links carry more signal than concepts. We also report where it does not help: on concepts a plain average of the same evidence is slightly better (\cref{sec:ablation,sec:types}).
  \item \textbf{A link between emotional content and epistemic reliability}: emotionally loaded concepts receive lower reliability, and confidence detects hallucination less well on them (\cref{sec:emotion}).
\end{enumerate}

\section{Related work}
\label{sec:related}

\paragraph{Hallucination detection} Surveys organise the field by where the signal comes from \citep{ji2023survey}. Internal detectors use sequence likelihood and entropy \citep{malinin2021uncertainty}, self-evaluation \citep{kadavath2022language} or probes on hidden states \citep{azaria-mitchell-2023-internal,chen2024inside}. SelfCheckGPT \citep{manakul-etal-2023-selfcheckgpt} scores sentences by their consistency with resamples; its NLI variant is a strong zero-resource baseline and the one we compare against. Semantic entropy clusters resamples by meaning \citep{kuhn2023semantic,farquhar2024detecting}. Grounding detectors score entailment between source and output, from FactCC \citep{kryscinski-etal-2020-evaluating} to SummaC \citep{laban-etal-2022-summac}, TRUE \citep{honovich-etal-2022-true} and AlignScore \citep{zha-etal-2023-alignscore}. FActScore \citep{min-etal-2023-factscore} decomposes long answers into atomic facts. \zfmn{} differs in its unit, concepts and the links between them, and in how it combines detectors: as evidence with an explicit amount, not as scores to average.

\paragraph{Benchmarks} RAGTruth \citep{niu-etal-2024-ragtruth} provides word-level spans with four hallucination types (evident or subtle, conflict or baseless) for answers of six models on question answering, summarisation and data-to-text generation. Mu-SHROOM \citep{vazquez-etal-2025-semeval} offers multilingual span labels. To our knowledge no benchmark labels the concepts and relations of an answer, which is what \bench{} adds.

\paragraph{Forma mentis networks} TFMNs \citep{stella2020text} model text as networks of lemmatised concepts linked by syntactic proximity, and EmoAtlas \citep{semeraro2025emoatlas} attaches emotional and valence profiles to their nodes using lexicons such as the NRC Emotion Lexicon \citep{mohammad2013crowdsourcing}. They have been used to study how people and models frame concepts. We use the same construction but attach epistemic confidence to every node and edge.

\paragraph{Uncertainty formalisms} Z-numbers \citep{zadeh2011note} pair a fuzzy restriction $A$ on a value with a fuzzy reliability $B$; \citet{kang2012method} convert them into classical fuzzy numbers, and later work develops their arithmetic \citep{yager2012z,aliev2015arithmetic}. Subjective logic \citep{josang2016subjective} represents an opinion as belief, disbelief and uncertainty derived from evidence counts, with fusion operators for combining sources. Evidential deep learning \citep{sensoy2018evidential} uses such opinions as the output of a trained classifier. We use them instead as a training-free calculus for fusing heterogeneous detector outputs on a graph.

\section{Method}
\label{sec:method}

\subsection{The answer as a network}
We parse the answer with spaCy \citep{honnibal2020spacy} and keep lemmatised content words as nodes. Two nodes are linked when they lie within three dependency hops in the same sentence, following the TFMN construction \citep{stella2020text}. We score the answer with the model that generated it (teacher forcing) and align each token's log-probability to the word it belongs to. A node's confidence $c_v$ is the geometric mean token probability of its word, and an edge's confidence $c_e$ is the minimum word confidence on its dependency path, so a link is only as confident as its weakest word. The same construction is applied to $N=10$ resamples of the answer.

\subsection{Subjective-logic opinions}
An opinion $\omega=(b,d,u)$ with $b+d+u=1$ expresses belief, disbelief and uncertainty about a proposition, here ``this element is supported''. It is built from $r$ units of positive and $s$ units of negative evidence with a prior weight $W$:
\begin{equation}
  b=\frac{r}{\lambda},\quad d=\frac{s}{\lambda},\quad u=\frac{W}{\lambda},\quad \lambda=r+s+W.
  \label{eq:evidence}
\end{equation}
The projected probability is $P=b+a\,u$ with base rate $a=\tfrac12$. With little evidence $u$ is large and $P$ stays near $a$; with much evidence $u$ shrinks and $P$ follows the evidence ratio. We use $W=2$, which makes the opinion equivalent to a Beta posterior under a uniform prior \citep{josang2016subjective}.

\subsection{Three evidence sources}
\paragraph{Model opinion $\omega^{\mathrm{M}}$} Over the $T=N+1$ texts (the answer and its resamples), the element contributes positive evidence equal to its confidence in each text where it occurs and zero where it does not: $r=\sum_{t}c_t$, $s=T-r$. An element that the model writes confidently and repeats across resamples gains belief.

\paragraph{Source opinion $\omega^{\mathrm{S}}$} Each answer sentence is checked against the source with DeBERTa-large-MNLI \citep{he2021deberta,williams-etal-2018-broad}. Following SummaC-ZS \citep{laban-etal-2022-summac}, the source is split into chunks of at most 300 words and the chunk with the highest entailment minus contradiction is kept. Its entailment and contradiction probabilities $e$ and $c$ become evidence worth $k$ observations, $r=k\,e$ and $s=k\,c$, with $k=10$. The neutral mass $1-e-c$ adds no evidence, so a neutral verdict raises $u$ instead of moving $P$. An element that occurs in several sentences takes the evidence of its most contradicted occurrence.

\paragraph{Consistency opinion $\omega^{\mathrm{C}}$} The SelfCheckGPT-NLI score $\sigma\in[0,1]$ of the element's sentence (the mean over the $N$ resamples of $c/(c+e)$ with each resample as premise, \citealp{manakul-etal-2023-selfcheckgpt}) becomes $r=N(1-\sigma)$ and $s=N\sigma$.

\subsection{Cumulative fusion}
Two opinions on the same proposition from independent sources are combined by cumulative fusion \citep{josang2016subjective}:
\begin{equation}
  \begin{gathered}
  b_{1\oplus2}=\frac{b_1u_2+b_2u_1}{\kappa},\qquad
  u_{1\oplus2}=\frac{u_1u_2}{\kappa},\\
  \kappa=u_1+u_2-u_1u_2,
  \end{gathered}
  \label{eq:fusion}
\end{equation}
with $d_{1\oplus2}=1-b_{1\oplus2}-u_{1\oplus2}$. Each source's belief is weighted by the \emph{other} source's uncertainty, so a source with much evidence (low $u$) dominates one with little. This is what separates fusion from an average with fixed weights: when the entailment model is neutral, $\omega^{\mathrm{S}}$ has high $u$ and the fused opinion defers to the model and consistency evidence; when it is confident, it dominates. The full opinion is $\omega=(\omega^{\mathrm{M}}\oplus\omega^{\mathrm{S}})\oplus\omega^{\mathrm{C}}$ (the operator is commutative and associative).

\subsection{Z-number and scores}
The element's Z-number is $Z=(A,B)$. $A$ is a triangular fuzzy number over the element's confidence in the eleven texts (minimum, median, maximum), and $B$ is the reliability, a triangular fuzzy number with support $[b,\,b+u]$ and peak $P$ of the fused opinion. Kang's conversion \citep{kang2012method} turns $Z$ into one fuzzy number when a crisp value is needed. The hallucination risk of an element is $1-P$, and the risk of an answer is the maximum over its links. No step uses labels; the only free parameters ($W=2$, $k=10$, the three-hop window and $N=10$) were fixed before evaluation, and \cref{sec:sensitivity} shows the result is stable around them.

\section{The \bench{} benchmark}
\label{sec:bench}

\paragraph{Construction} We take all 450 RAGTruth test answers of each of Mistral-7B-Instruct \citep{jiang2023mistral}, Llama-2-7B-chat and Llama-2-13B-chat \citep{touvron2023llama}, 150 per task, and build their TFMNs. An element is labelled hallucinated when any of its occurrences overlaps a human-annotated hallucination span, and it inherits that span's type. A link is covered when either of its words lies in the span. Each element row also holds its degree, its EmoAtlas valence and emotion features and the score of every method in this paper; each answer row holds the answer label and answer-level scores.

\paragraph{Statistics} \cref{tab:bench} summarises the release. Half of the answers (684 of 1,350) contain a hallucination, but only 7.4\% of links and 7.6\% of concepts are hallucinated, so element-level detection is strongly imbalanced, and we report AUROC as the main metric with precision--recall curves alongside. Data-to-text answers are hallucinated most often at the answer level; summaries have the lowest element rates.

\begin{table}[!htbp]
  \centering
  \caption{\textbf{\bench{} statistics.} Nine datasets (model $\times$ task), 71,674 concepts and 201,782 links in total. D2T: data-to-text. ``hal.'' counts answers with at least one annotated hallucination and the share of hallucinated elements.}
  \label{tab:bench}
  \small
  \setlength{\tabcolsep}{3.4pt}
  \begin{adjustbox}{max width=\textwidth}
\begin{tabular}{llrrrrr}
    \toprule
    & & \multicolumn{2}{c}{Answers} & \multicolumn{2}{c}{Links} & Concepts \\
    \cmidrule(lr){3-4}\cmidrule(lr){5-6}\cmidrule(lr){7-7}
    Model & Task & all & hal. & all & hal.\,\% & hal.\,\% \\
    \midrule
    \multirow{3}{*}{Llama-2-7B}  & QA      & 150 & 52  & 24,126 & 8.6 & 10.1 \\
                                 & Summary & 150 & 51  & 19,645 & 5.4 & 5.2 \\
                                 & D2T     & 150 & 123 & 26,751 & 8.6 & 7.7 \\
    \midrule
    \multirow{3}{*}{Mistral-7B}  & QA      & 150 & 31  & 13,935 & 6.6 & 7.3 \\
                                 & Summary & 150 & 86  & 21,104 & 8.4 & 8.6 \\
                                 & D2T     & 150 & 134 & 24,065 & 7.8 & 8.7 \\
    \midrule
    \multirow{3}{*}{Llama-2-13B} & QA      & 150 & 36  & 20,560 & 7.4 & 7.9 \\
                                 & Summary & 150 & 33  & 18,332 & 3.0 & 2.7 \\
                                 & D2T     & 150 & 138 & 33,264 & 8.8 & 8.9 \\
    \midrule
    \multicolumn{2}{l}{Total} & 1,350 & 684 & 201,782 & 7.4 & 7.6 \\
    \bottomrule
  \end{tabular}
\end{adjustbox}
\end{table}

\begin{figure}[!htbp]
  \centering
  \includegraphics[width=0.62\textwidth]{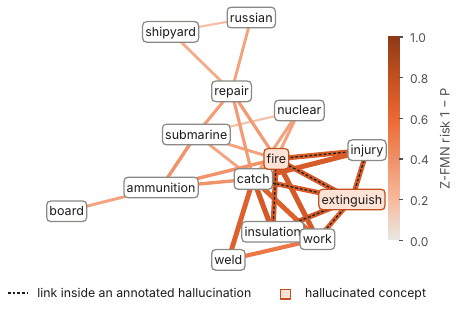}
  \caption{\textbf{One answer as a network.} A Llama-2-7B news summary (14 concepts, 30 links). Edge colour and width show \zfmn{} risk; dashed links lie inside the human-annotated hallucination. The risky links cluster around the invented event (annotated as evident baseless information) (\emph{fire}, \emph{extinguish}, \emph{injury}), while the correct background (\emph{submarine}, \emph{shipyard}, \emph{repair}) stays low. Within this answer, risk ranks hallucinated links with AUROC 0.84.}
  \label{fig:example}
\end{figure}

\paragraph{Licence and release} RAGTruth is released under the MIT licence. \bench{} contains labels, features, scores and RAGTruth response identifiers; source passages and answers are not redistributed and join back to RAGTruth by identifier. \cref{fig:example} shows one answer as it appears in the benchmark.

\section{Experimental setup}
\label{sec:setup}

\paragraph{Models and generation} Each answer is scored by the model that wrote it, loaded in 4-bit precision \citep{dettmers2023qlora}, and $N=10$ resamples are drawn from the original prompt at temperature 0.7. All runs used Kaggle's free tier (two T4 GPUs), roughly 60 GPU-hours in total.

\paragraph{Baselines} (i) Mean token log-probability of the element. (ii) SelfCheckGPT-NLI \citep{manakul-etal-2023-selfcheckgpt} with the same resamples and the same DeBERTa-large-MNLI model, applied to the element's sentence. (iii) Source NLI: $1-e$ of the element's sentence against the source, using the same chunking. (iv) The model opinion alone. (v) Training-free combinations of the three risks $1-P^{\mathrm{M}}$, $1-P^{\mathrm{S}}$ and $\sigma$: their mean and their maximum. (vi) Logistic stacking of the same six evidence features ($P^{\mathrm{M}}$, the context-conditioned $P$ and $u$, $e$, $c$, $\sigma$), trained with five-fold cross-validation grouped by answer. Stacking uses the labels that \zfmn{} never sees and serves as a supervised reference.

\paragraph{Metrics} AUROC for hallucinated links, concepts and answers; average precision via precision--recall curves; expected calibration error (ECE, 15 bins; \citealp{guo2017calibration}) and Brier score.

\paragraph{Statistics} A dataset is one model $\times$ task pair, nine in all. Confidence intervals (95\%) come from 1,000--2,000 bootstrap resamples of whole answers \citep{efron1994introduction}, so all elements of an answer stay together. Every comparison is paired on the same resample. Pooled results average the per-dataset AUROCs (or their differences), with the bootstrap stratified by dataset. We call a difference significant when its interval excludes zero.

\section{Results}
\label{sec:results}

\subsection{Main comparison}
\label{sec:main}

\begin{figure}[!htbp]
  \centering
  \includegraphics[width=\textwidth]{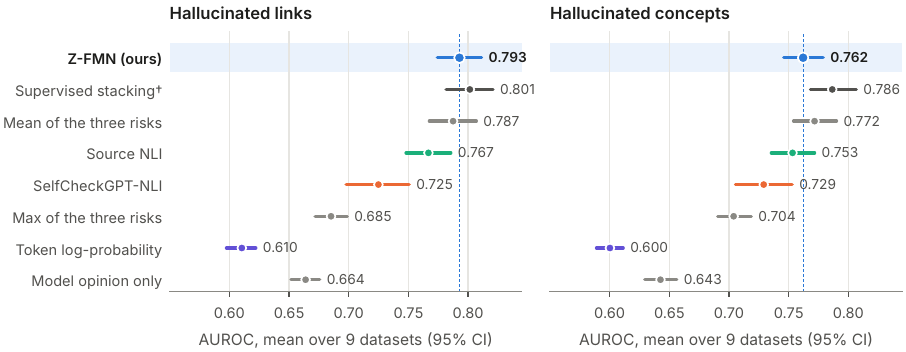}
  \caption{\textbf{Main result.} AUROC for hallucinated links (left) and concepts (right), averaged over the nine model $\times$ task datasets, with 95\% bootstrap intervals over answers. \zfmn{} is the best training-free detector on links; on concepts the plain average of the same three risks is slightly higher. $^\dagger$Supervised stacking is trained on the labels with cross-validation and is shown as a reference, not a competitor.}
  \label{fig:main}
\end{figure}

\begin{table}[!htbp]
  \centering
  \caption{\textbf{Paired differences} between \zfmn{} and each method, pooled over datasets (95\% CI). Bold: \zfmn{} significantly better. Answer level: risk of the riskiest link against each method's answer score (\cref{sec:answer}). $^\dagger$Supervised.}
  \label{tab:main}
  \small
  \setlength{\tabcolsep}{2.6pt}
  \begin{adjustbox}{max width=\textwidth}
\begin{tabular}{lccc}
    \toprule
    & \multicolumn{3}{c}{\zfmn{} $-$ method (AUROC)} \\
    \cmidrule(lr){2-4}
    Method & Links & Concepts & Answers \\
    \midrule
    SelfCheckGPT-NLI & \best{+0.068}\ci{.053}{.083} & \best{+0.033}\ci{.018}{.047} & \best{+0.049}\ci{.020}{.075} \\
    Source NLI & \best{+0.026}\ci{.013}{.040} & +0.009\ci{$-$.006}{.024} & +0.007\ci{$-$.037}{.052} \\
    Token log-prob & \best{+0.182}\ci{.162}{.202} & \best{+0.162}\ci{.143}{.182} & \best{+0.191}\ci{.129}{.249} \\
    Model opinion only & \best{+0.129}\ci{.112}{.145} & \best{+0.120}\ci{.105}{.135} & -- \\
    Mean of risks & \best{+0.005}\ci{.002}{.009} & $-$0.010\ci{$-$.013}{$-$.006} & -- \\
    Max of risks & \best{+0.108}\ci{.092}{.123} & \best{+0.058}\ci{.050}{.066} & -- \\
    \midrule
    Stacking$^\dagger$ & $-$0.009\ci{$-$.017}{$-$.001} & $-$0.024\ci{$-$.034}{$-$.016} & -- \\
    \bottomrule
  \end{tabular}
\end{adjustbox}
\end{table}

\cref{fig:main} and \cref{tab:main} give the headline result. On links, \zfmn{} reaches AUROC 0.793 \ci{0.774}{0.811}, against 0.767 for source NLI, 0.725 for SelfCheckGPT-NLI and 0.610 for token log-probability. Every training-free competitor is significantly worse. The closest is the plain mean of the same three risks (0.787), which uses exactly the same evidence and trails by 0.005. On concepts, \zfmn{} (0.762) is significantly ahead of SelfCheckGPT and token probability and level with source NLI, but the plain mean is ahead of it by 0.010; we return to this in \cref{sec:ablation}.

The supervised stacking model, which sees the labels, is ahead by 0.009 on links and 0.024 on concepts. \zfmn{} therefore recovers most of what a trained combiner extracts from these features without any labelled data, which matters for a detector meant to be deployed on new models and domains.

\begin{figure}[!htbp]
  \centering
  \includegraphics[width=0.72\textwidth]{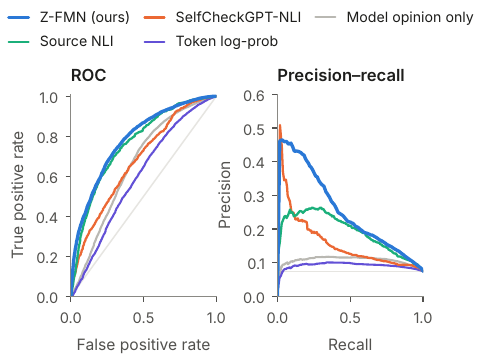}
  \caption{\textbf{ROC and precision--recall curves} for hallucinated links, all datasets pooled (201,782 links, 7.4\% hallucinated; the grey line in the right panel is this base rate). \zfmn{} dominates over most of the ROC curve and keeps precision near 0.45 up to 10\% recall, where SelfCheckGPT drops quickly.}
  \label{fig:roc}
\end{figure}

The curves in \cref{fig:roc} show where the gain comes from. At low false-positive rates, the regime that matters for flagging claims to a reader, \zfmn{} keeps precision near 0.45 over the first 10\% of recall, about six times the base rate, while SelfCheckGPT's precision falls below 0.3 almost at once and source NLI stays near 0.25.

\subsection{Robustness across models, tasks and scale}
\label{sec:robust}

\begin{figure}[!htbp]
  \centering
  \includegraphics[width=0.72\textwidth]{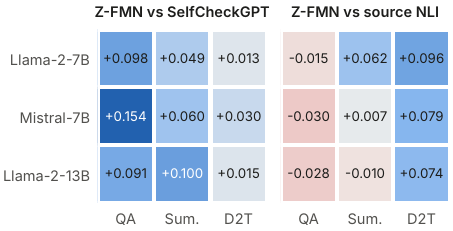}
  \caption{\textbf{Gain per dataset.} Link AUROC of \zfmn{} minus SelfCheckGPT-NLI (left) and minus source NLI (right). Blue: \zfmn{} better. \zfmn{} beats SelfCheckGPT on all nine datasets. Against source NLI it gains most on data-to-text and loses slightly on QA, where entailment against a short passage is already strong.}
  \label{fig:heat}
\end{figure}

\cref{fig:heat} breaks the link result down by dataset. \zfmn{} beats SelfCheckGPT on all nine, by 0.013 to 0.154. Against source NLI the picture is task-dependent: \zfmn{} gains 0.074--0.096 on data-to-text, where the source is structured JSON and entailment is least reliable, and trails by 0.015--0.030 on QA, where a short passage makes entailment strong on its own (and by 0.010 on Llama-2-13B summaries). This is the behaviour fusion is designed for: when the source opinion is uncertain, the other sources take over.

The gain holds across scale. Pooled over tasks, the margin over SelfCheckGPT on links is $+0.054$ \ci{.029}{.081} for Llama-2-7B, $+0.081$ \ci{.054}{.111} for Mistral-7B and $+0.069$ \ci{.040}{.096} for Llama-2-13B. Moving from 7B to 13B within the Llama-2 family leaves \zfmn{}'s AUROC unchanged (0.768 and 0.767), while SelfCheckGPT drops from 0.714 to 0.698.

\subsection{What the Z-number contributes}
\label{sec:ablation}

\begin{figure}[!htbp]
  \centering
  \includegraphics[width=0.72\textwidth]{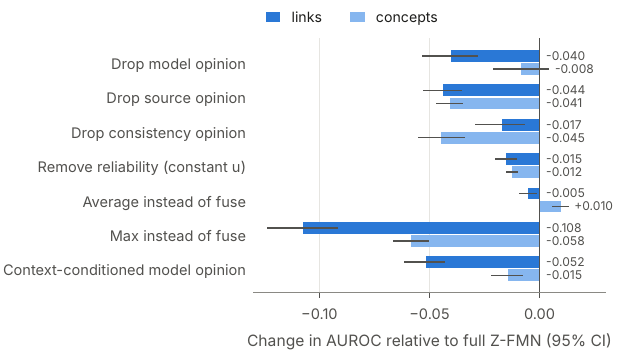}
  \caption{\textbf{Ablations.} Change in AUROC when one component of \zfmn{} is removed or replaced (pooled over datasets, 95\% CI). Every source helps on links, and making each source's uncertainty constant, which keeps its probability but removes its reliability, costs 0.015. Replacing fusion by a plain average costs little on links and helps slightly on concepts; replacing it by the maximum costs 0.108.}
  \label{fig:ablation}
\end{figure}

\cref{fig:ablation} removes one component at a time.

\paragraph{Every source contributes} Dropping the model opinion costs 0.040 on links, dropping the source opinion 0.044 and dropping the consistency opinion 0.017. On concepts, the source and consistency opinions matter most (0.041 and 0.045), while the model opinion adds little.

\paragraph{Fusion versus averaging} Replacing cumulative fusion with the plain mean of the three risks costs 0.005 \ci{.002}{.009} on links and gains 0.010 on concepts; replacing it with their maximum costs 0.108 and 0.058. In AUROC, then, a well-chosen average is a strong rival. Fusion has two advantages that AUROC does not show: it needs no choice of weights, since each source is weighted by its own evidence, and it returns an opinion $(b,d,u)$ rather than a score, which is what makes the type analysis of \cref{sec:types} possible.

\paragraph{Reliability matters} The sharpest test of the Z-number's reliability part keeps every source's projected probability but replaces its uncertainty with a constant, so the fusion still runs but cannot tell a confident source from an unsure one. This costs 0.015 \ci{.011}{.020} on links and 0.012 \ci{.010}{.015} on concepts, significant on both. Fusion without reliability (0.778) is in fact worse than the plain average (0.787); it is the evidence amounts carried by the B component that lift it above.

\paragraph{A richer model opinion did not help} A ``context-conditioned'' variant counts a resample that contains an element's neighbours (or, for a link, its endpoints) but not the element itself as negative evidence. Its uncertainty varies per element, yet it lowers AUROC by 0.052. The simpler opinion is better.

\subsection{Which hallucinations does it catch?}
\label{sec:types}

\begin{figure}[!htbp]
  \centering
  \includegraphics[width=0.72\textwidth]{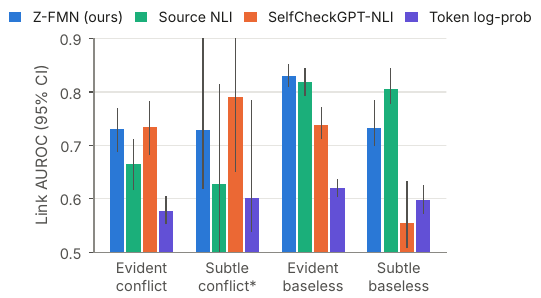}
  \caption{\textbf{By hallucination type.} Link AUROC against correct links for each RAGTruth type. \zfmn{} is best or level on evident errors, ties SelfCheckGPT on evident conflicts, and is ahead of SelfCheckGPT by 0.18 on subtle baseless content, where source NLI is strongest. *Subtle conflicts occur in only three datasets, hence the wide intervals.}
  \label{fig:types}
\end{figure}

\paragraph{Invented versus contradicted content} \cref{fig:types} splits the hallucinated links by type. On invented (baseless) content \zfmn{} is far ahead of SelfCheckGPT: $+0.105$ \ci{.088}{.124} pooled, and $+0.18$ on subtle baseless links, which resamples rarely contradict because the model invents them consistently. On contradictions the two are level ($-0.007$ \ci{$-$.033}{.026}), and both beat source NLI clearly (0.731 and 0.735 against 0.666 on evident conflicts). Contradictions are where the consistency opinion earns its place: adding it to the model and source opinions raises conflict AUROC by 0.064 \ci{.040}{.086}. Source NLI remains the best single detector of subtle baseless content (0.806), a case we return to in \cref{sec:discussion}.

\begin{figure}[!htbp]
  \centering
  \includegraphics[width=\textwidth]{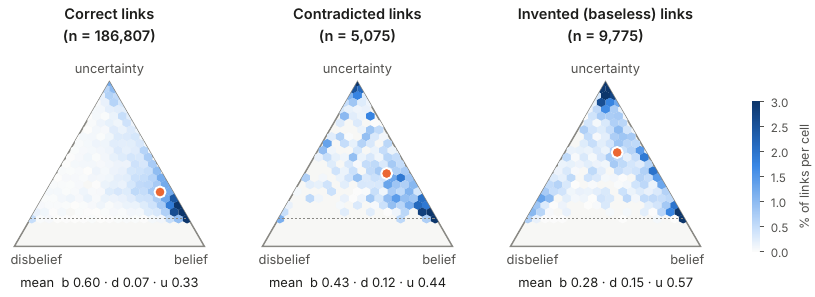}
  \caption{\textbf{The source opinion of every link in the opinion simplex.} Each triangle shows where links fall between belief, disbelief and uncertainty (colour: share of links per cell; the dotted line is the uncertainty floor set by $k=10$). Correct links concentrate near belief. Contradicted links move away from it, and invented links move toward \emph{uncertainty}: the entailment model is neutral about content the source never mentions. A scalar probability places both kinds in the same middle range; the simplex keeps them apart.}
  \label{fig:simplex}
\end{figure}

\paragraph{The uncertainty dimension tells the two apart} \cref{fig:simplex} shows why a three-valued opinion helps. Correct links sit near belief (mean $b=0.60$, $u=0.33$). Contradicted links lose belief ($b=0.43$, $u=0.44$), and invented links are dominated by uncertainty ($b=0.28$, $u=0.57$), because the entailment model returns ``neutral'' for content the source never mentions. Among hallucinated links, low fused uncertainty identifies contradictions against invented content with AUROC 0.704 (concepts 0.691), while the scalar risk $1-P$ reaches only 0.602 in its better orientation and NLI contradiction 0.582. The separation is inherited from the entailment model's neutral mass (its AUROC is identical), so the claim is not that fusion creates this information but that the Z-number \emph{keeps} it while still producing the best training-free link score; collapsing to $P$ loses about half of its separation above chance. A reader of \zfmn{}'s output can therefore be told not only that a link is doubtful but whether it contradicts the source or is unsupported by it.

\subsection{Links carry more signal than concepts}
\label{sec:links}
With the same method, links are easier to verify than concepts: \zfmn{} scores 0.031 \ci{.023}{.039} higher on links, and the model-plus-source fusion 0.058 higher. The gap is absent for SelfCheckGPT ($-0.004$, n.s.), whose score is the same for every element of a sentence. The advantage is therefore not a property of the labels but of the relational evidence: a link inherits the doubt of its weakest word and is checked as part of a specific claim.

\subsection{Answer-level detection}
\label{sec:answer}
Scoring an answer by its riskiest link gives AUROC 0.672 \ci{.629}{.710}, ahead of SelfCheckGPT (0.623; $+0.049$) and token log-probability (0.481), and level with source NLI (0.665; $+0.007$, n.s.). \zfmn{}'s strength is the element level, where it says \emph{where} an answer goes wrong; at the answer level it matches the best grounding baseline. Per-dataset values are in \ref{app:tables}.

\subsection{Sensitivity to the evidence weights}
\label{sec:sensitivity}

\begin{figure}[!htbp]
  \centering
  \begin{subfigure}[b]{0.305\textwidth}
    \centering
    \includegraphics[width=\textwidth]{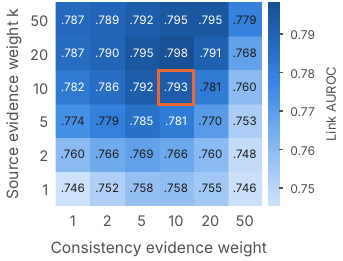}
    \caption{Evidence weights}
    \label{fig:sens}
  \end{subfigure}\hfill
  \begin{subfigure}[b]{0.25\textwidth}
    \centering
    \includegraphics[width=\textwidth]{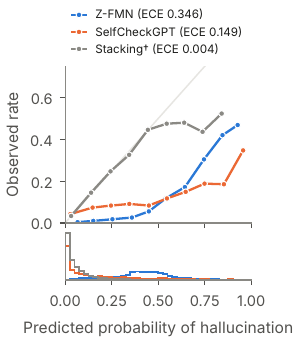}
    \caption{Calibration}
    \label{fig:calib}
  \end{subfigure}\hfill
  \begin{subfigure}[b]{0.43\textwidth}
    \centering
    \includegraphics[width=\textwidth]{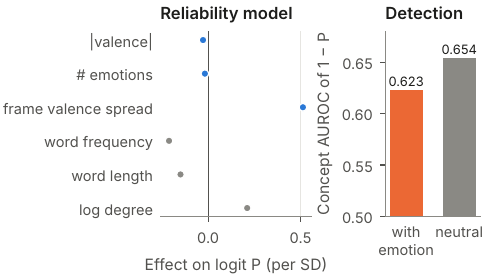}
    \caption{Emotion and reliability}
    \label{fig:emotion}
  \end{subfigure}
  \caption{\textbf{(a) Sensitivity:} link AUROC as a function of how many observations the source ($k$) and consistency opinions are worth; the box marks the setting used throughout ($k=10$, $N=10$), fixed before evaluation. \textbf{(b) Calibration:} reliability diagram on links with the distribution of predictions below. \zfmn{} ranks well but its probabilities are too high; the supervised stacker is calibrated by construction. \textbf{(c) Emotion} (71,674 concepts): left, a linear mixed model of the model opinion's logit $P$ with a random intercept per answer (coefficients per SD; 95\% CIs are narrower than the markers; emotional predictors in blue); right, the model opinion detects hallucinated concepts less well when they carry emotion.}
  \label{fig:sens_calib}
\end{figure}

The pre-chosen weights are not tuned on the test data, so we report how much they matter (\cref{fig:sens}). Across a 36-cell grid from 1 to 50 observations per source, link AUROC ranges from 0.746 to 0.798; the chosen setting gives 0.793, within 0.005 of the best cell, and every setting with $k\ge5$ and a consistency weight up to 20 scores at least 0.770, still ahead of source NLI and SelfCheckGPT. Giving the consistency opinion more weight than about 20 observations hurts, probably because it then overrides the source opinion on baseless content, where resamples agree with the invention.

\subsection{Calibration}
\label{sec:calibration}
\zfmn{} is a good ranker but a poorly calibrated probability (\cref{fig:calib}). Its link ECE, averaged over datasets, is 0.349, against 0.172 for SelfCheckGPT and 0.011 for the supervised stacker (the figure pools all links, giving similar values). The cause is the base rate: cumulative fusion starts from $a=\tfrac12$, while only 7\% of links are hallucinated, so moderate risks are overstated. Where calibrated probabilities are needed, a monotone recalibration on a small labelled set fixes this without changing the ranking.

\subsection{Emotional content and reliability}
\label{sec:emotion}

EmoAtlas gives every concept a valence and a set of NRC emotions, so we can ask whether the model's reliability depends on emotional content (\cref{fig:emotion}). Controlling for word frequency, length, degree, task and model, concepts with stronger valence receive lower reliability ($-0.028$ SD \ci{$-$.039}{$-$.016}, $p<10^{-5}$), as do concepts with more emotions ($-0.018$, $p=0.003$). The strongest predictor is the valence spread of a concept's neighbours ($+0.51$ SD): concepts framed by a mix of positive and negative words are held more confidently, possibly because they sit in richer, more repeated contexts. Confidence also detects hallucination less well on emotional concepts: the interaction between risk and emotional content has odds ratio 0.895 \ci{0.830}{0.964} ($p=0.003$), and the AUROC is 0.623 on emotional against 0.654 on neutral concepts. Emotionally loaded claims are thus both less trusted by the model and harder to verify from its confidence alone.

\section{Engineering application: claim-level fault monitoring}
\label{sec:deploy}

\paragraph{The monitor in a deployed system} \zfmn{} sits after the generator of a retrieval-augmented system. For each answer it builds the concept network, collects the three opinions, fuses them and returns the answer's links ranked by risk $1-P$, each with its opinion $(b,d,u)$. A reviewer, or a rule in the downstream pipeline, checks the riskiest links first. The opinion also says what kind of check is needed: a link with high disbelief is contradicted by the source and can be compared with the passage directly, while a link with high uncertainty is not mentioned by the source at all and needs an outside check (\cref{sec:types}). No step uses labels, so the monitor can be attached to a new model or domain as it is.

\begin{figure}[!htbp]
  \centering
  \includegraphics[width=\textwidth]{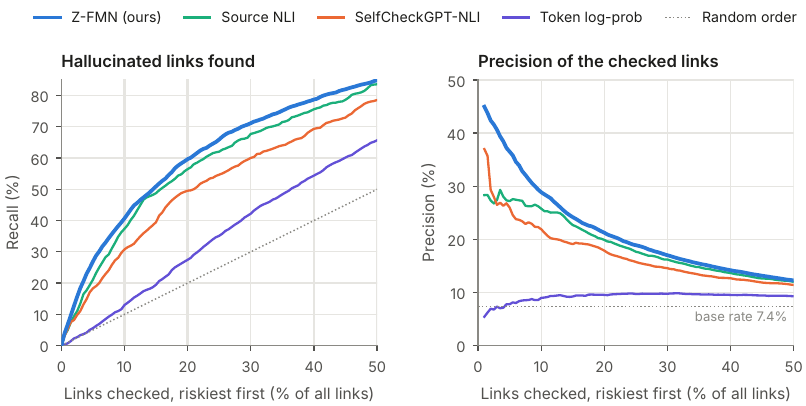}
  \caption{\textbf{Fault monitoring under a review budget.} A reviewer checks links in order of risk, riskiest first. Left: share of all hallucinated links found; right: share of the checked links that are hallucinated. Mean over the nine model $\times$ task datasets. The grey lines on the left mark review budgets of 10\% and 20\% of links.}
  \label{fig:budget}
\end{figure}

\begin{table}[!htbp]
  \centering
  \caption{\textbf{Review-budget results.} Recall: share of hallucinated links found when the riskiest $q\%$ of links are checked. Precision: share of the checked links that are hallucinated (base rate 7.4\%). Mean over the nine datasets; bold marks the best method in each column.}
  \label{tab:budget}
  \small
  \begin{tabular}{lcccccc}
    \toprule
    & \multicolumn{4}{c}{Recall (\%) at budget $q$} & \multicolumn{2}{c}{Precision (\%) at $q$} \\
    \cmidrule(lr){2-5}\cmidrule(lr){6-7}
    Method & 5\% & 10\% & 20\% & 30\% & 5\% & 10\% \\
    \midrule
    \ours \zfmn{} (ours) & \best{26.9} & \best{40.7} & \best{59.7} & \best{71.2} & \best{36.5} & \best{28.9} \\
    Source NLI & 21.1 & 37.4 & 56.5 & 67.7 & 27.2 & 25.8 \\
    SelfCheckGPT-NLI & 18.3 & 30.8 & 49.6 & 60.1 & 26.0 & 21.9 \\
    Token log-prob & 5.6 & 13.0 & 27.5 & 42.2 & 7.8 & 9.0 \\
    Random order & 5.0 & 10.0 & 20.0 & 30.0 & 7.4 & 7.4 \\
    \bottomrule
  \end{tabular}
\end{table}

\paragraph{How much a reviewer finds} \cref{fig:budget} and \cref{tab:budget} give the share of hallucinated links a reviewer finds when checking only the riskiest fraction of links. Checking 10\% of links finds 40.7\% of the hallucinated ones with \zfmn{}, against 37.4\% with source NLI, 30.8\% with SelfCheckGPT-NLI and 10\% in random order. At a 20\% budget the shares are 59.7\%, 56.5\% and 49.6\%. The gain over SelfCheckGPT is 9.9 percentage points at 10\% \ci{7.1}{14.3} and 10.1 points at 20\% \ci{6.9}{13.2}; the gain over source NLI is smaller, 3.3 \ci{0.0}{6.7} and 3.2 points \ci{0.1}{6.7}, with intervals from 1,000 bootstrap resamples of answers. Among the riskiest 5\% of links, 36.5\% are hallucinated, almost five times the base rate, so a reviewer who starts at the top of the list spends most of the early effort on real faults.

\paragraph{Cost per answer} The monitor needs the generator's token probabilities, ten resampled answers and sentence-level entailment. On the free hardware used here (Kaggle's free tier with two NVIDIA T4 GPUs, 4-bit models), generating and scoring an answer with its resamples took about 65, 103 and 125 seconds for Llama-2-7B, Mistral-7B and Llama-2-13B, and entailment about 9 seconds; network construction, fusion and scoring take milliseconds on a CPU (\ref{app:impl}). Resampling dominates the cost, so a system that already draws resamples for SelfCheckGPT adds only the entailment step. Where latency matters, the monitor can run asynchronously on logged answers, or with fewer resamples at some loss in the consistency opinion.

\paragraph{Choosing the operating point} The fused probability is a good ranking but not a calibrated probability (\cref{sec:calibration}), so a fixed threshold on $P$ does not transfer between systems. Two practical choices avoid this: fix the review budget and check the top $q\%$ of links per batch, as above, or fit a monotone recalibration on a small labelled sample from the target system, which changes the threshold but not the ranking.

\section{Discussion}
\label{sec:discussion}

\paragraph{What the Z-number buys} Two findings support representing confidence as a Z-number rather than a probability. Removing reliability while keeping the probabilities costs a significant 0.015, so the evidence amounts matter; and the uncertainty dimension separates contradicted from invented content, information a probability discards. The detection gain over the best training-free combination is small (0.005 on links, and none on concepts), so the case for the Z-number rests less on AUROC than on what it reports: a weight-free combination whose output says how much evidence there is and of which kind. The value is not in any single source, each of which is a known detector, but in combining them with an explicit account of how much each one knows.

\paragraph{What it does not buy} The model opinion's uncertainty is the same for every element, because it depends only on the number of texts ($u=W/(W+T)=0.154$). Its projected probability is therefore equivalent in ranking to the plain resample frequency (difference $-0.003$, n.s.), and all variation in reliability comes from the source and consistency opinions. Making the model opinion's evidence weight depend on the element, for instance on how often it appears in pretraining data, is the most direct way to strengthen this part. Source NLI alone remains better on QA and on subtle baseless content, which suggests weighting the source opinion by task, something a small labelled set could estimate.

\paragraph{Why links} Hallucinations in RAGTruth are mostly wrong combinations of right words. A link carries the confidence of its weakest word and the verdict on the claim it belongs to, which is why it is the better unit for every detector that is not sentence-constant. For users, a flagged link (``\emph{fire}--\emph{extinguish}'') points at the claim to check more precisely than a flagged sentence.

\section{Conclusion}
\label{sec:conclusion}
We asked whether confidence about an answer is better placed on its concepts and links than on its tokens or sentences, and whether it should be represented as evidence with an explicit reliability rather than as one probability. \zfmn{} does both. Without any training, it detects hallucinated links better than SelfCheckGPT and source entailment on all three models, at least matches simple combinations of the same evidence, comes close to a supervised combiner, and tells contradictions from invented content. We release \bench{} so that future detectors can be compared on where in an answer they place doubt, not only on whether they flag it.

\section{Limitations and future work}
\begin{itemize}
  \item \textbf{One benchmark, one language.} All results are on RAGTruth in English with three open 7--13B models. A second span-labelled benchmark, such as Mu-SHROOM \citep{vazquez-etal-2025-semeval}, is needed before claiming generality across languages and domains.
  \item \textbf{Access to the generator.} \zfmn{} reads token probabilities and draws resamples. Closed models would need a proxy scorer, which we have not tested.
  \item \textbf{Entailment errors.} The NLI model was trained on natural sentence pairs and is weaker on structured data-to-text sources; the source opinion inherits its errors.
  \item \textbf{Independence assumption.} Cumulative fusion assumes independent evidence. The model and consistency opinions share the same resamples, so the fused uncertainty is probably too small.
  \item \textbf{Calibration.} Raw $P$ is not a calibrated probability (\cref{sec:calibration}).
  \item \textbf{Label projection.} Element labels are projected from word spans: a link is hallucinated when either word lies in a span, which can mark a correct relation that touches a wrong word.
  \item \textbf{Cost.} Each answer needs ten resamples and sentence-level NLI, about 1--2 GPU-minutes on a T4, far more than single-pass detectors (\cref{sec:deploy}).
\end{itemize}
Future work will test the monitor on newer and closed models through a proxy scorer, on a second benchmark, and with a model opinion whose evidence weight depends on the element.

\FloatBarrier
\section*{Declaration of competing interest}
The author(s) declare that they have no known competing financial interests or personal relationships that could have appeared to influence the work reported in this paper.

\section*{Funding}
This research did not receive any specific grant from funding agencies in the public, commercial, or not-for-profit sectors.

\section*{Ethical statement}
This work uses only public secondary data. RAGTruth is released under the MIT licence, and its existing human annotations are used without collecting new ones. No human participants were involved and no personal data was collected. The detector is meant to flag unsupported claims for human review. It should not be used to certify an answer as true: a link with high $P$ can still be wrong, and lower reliability on emotionally loaded concepts means such claims deserve extra scrutiny rather than automatic trust or rejection.

\section*{Data availability}
RAGTruth is publicly available under the MIT licence. The \bench{} benchmark (element- and answer-level labels, features and the scores of every method) and the code that reproduces every table and figure are available at \ifanon[repository link removed for anonymised review]\else\url{https://github.com/SwetaKaman/z-fmn}\fi. The released files contain identifiers, not the source passages, which join back to RAGTruth.

\ifanon\else
\section*{CRediT authorship contribution statement}
\textbf{Sweta Kaman:} Conceptualization, Methodology, Software, Formal analysis, Investigation, Data curation, Visualization, Writing -- original draft, Writing -- review \& editing.

\section*{Acknowledgements}
Experiments used the free GPU tier of Kaggle.
\fi

\section*{Declaration of generative AI and AI-assisted technologies in the manuscript preparation process}
During the preparation of this work, the authors used Claude (Anthropic) in order to improve the grammatical structure, flow, and readability of the text. After using this tool, the authors reviewed, edited, and verified the final output as needed, and take full responsibility for the overall accuracy, integrity, and final content of the publication.

\bibliography{references}
\bibliographystyle{elsarticle-harv}

\appendix
\counterwithin{table}{section}
\renewcommand{\thetable}{\Alph{section}.\arabic{table}}

\section{Per-dataset results}
\label{app:tables}

\begin{table}[H]
  \centering
  \caption{\textbf{AUROC per dataset.} Link level (left) and answer level (right; answer risk is the riskiest link). Bold: best training-free method in the row. $^\dagger$Supervised reference, excluded from bold.}
  \label{tab:perdataset}
  \small
  \setlength{\tabcolsep}{4.5pt}
  \begin{adjustbox}{max width=\textwidth}
\begin{tabular}{llcccccccc}
    \toprule
    & & \multicolumn{5}{c}{Links} & \multicolumn{3}{c}{Answers} \\
    \cmidrule(lr){3-7}\cmidrule(lr){8-10}
    Model & Task & Log-prob & SelfCheck & NLI & \zfmn{} & Stacking$^\dagger$ & SelfCheck & NLI & \zfmn{} \\
    \midrule
    \multirow{3}{*}{Llama-2-7B} & QA & 0.613 & 0.677 & \best{0.790} & 0.775 & 0.802 & 0.560 & \best{0.562} & 0.538 \\
     & Summary & 0.591 & 0.722 & 0.709 & \best{0.771} & 0.765 & 0.664 & 0.751 & \best{0.763} \\
     & D2T & 0.534 & 0.743 & 0.661 & \best{0.756} & 0.759 & 0.622 & 0.481 & \best{0.650} \\
    \midrule
    \multirow{3}{*}{Mistral-7B} & QA & 0.684 & 0.729 & \best{0.913} & 0.883 & 0.914 & 0.677 & 0.697 & \best{0.815} \\
     & Summary & 0.658 & 0.788 & 0.842 & \best{0.849} & 0.857 & 0.604 & \best{0.746} & 0.722 \\
     & D2T & 0.613 & 0.769 & 0.720 & \best{0.799} & 0.804 & 0.608 & \best{0.699} & 0.640 \\
    \midrule
    \multirow{3}{*}{Llama-2-13B} & QA & 0.574 & 0.737 & \best{0.856} & 0.828 & 0.854 & 0.709 & 0.674 & \best{0.736} \\
     & Summary & 0.656 & 0.643 & \best{0.753} & 0.743 & 0.715 & 0.618 & \best{0.711} & 0.676 \\
     & D2T & 0.571 & 0.716 & 0.657 & \best{0.730} & 0.744 & 0.548 & \best{0.664} & 0.510 \\
    \midrule
    \multicolumn{2}{l}{Mean} & 0.610 & 0.725 & 0.767 & \best{0.793} & 0.801 & 0.623 & 0.665 & \best{0.672} \\
    \bottomrule
  \end{tabular}
\end{adjustbox}
\end{table}

\cref{tab:perdataset} lists link- and answer-level AUROC for every dataset. At the answer level, the clean-answer class is small for data-to-text (12--27 clean answers per model), so those AUROCs have wide intervals.

\section{Implementation details}
\label{app:impl}

\begin{table}[H]
  \centering
  \caption{Settings used for every result in the paper.}
  \label{tab:impl}
  \small
  \begin{adjustbox}{max width=\textwidth}
\begin{tabular}{ll}
    \toprule
    Setting & Value \\
    \midrule
    Generators & Mistral-7B-Instruct-v0.1, Llama-2-7B-chat, \\
               & Llama-2-13B-chat (4-bit NF4) \\
    Resamples $N$ & 10, temperature 0.7 \\
    Network & spaCy \texttt{en\_core\_web\_sm}, content \\
            & lemmas, links within 3 dependency hops \\
    NLI model & \texttt{microsoft/deberta-large-mnli}, fp32 \\
    Source chunks & $\le$300 words, best $e-c$ chunk \\
    Prior weight $W$, base rate $a$ & 2, 0.5 \\
    Source evidence $k$ & 10 \\
    Consistency evidence & $N=10$ \\
    Bootstrap & 1,000 (ablations), 2,000 (main tables) \\
    Hardware & Kaggle, 2$\times$T4, about 60 GPU-hours \\
    \bottomrule
  \end{tabular}
\end{adjustbox}
\end{table}

Generation with scoring took about 65, 103 and 125 seconds per answer for Llama-2-7B, Mistral-7B and Llama-2-13B, and NLI about 9 seconds per answer. All fusion, evaluation and statistics run on a CPU in minutes. The code reproduces every table and figure from the released element and answer files.

\section{Hallucination-type separation and calibration}
\label{app:types}

\begin{table}[H]
  \centering
  \caption{AUROC for telling contradicted from invented elements, among hallucinated elements only (0.5 = cannot tell them apart).}
  \label{tab:sep}
  \small
  \setlength{\tabcolsep}{4pt}
  \begin{adjustbox}{max width=\textwidth}
\begin{tabular}{lcc}
    \toprule
    Score (higher = contradiction) & Links & Concepts \\
    \midrule
    $-u$ (fused uncertainty) & \best{0.704} & \best{0.691} \\
    $-$NLI neutral mass & \best{0.704} & \best{0.691} \\
    Scalar $P$ (better orientation) & 0.602 & 0.619 \\
    NLI contradiction & 0.582 & 0.570 \\
    Disbelief $d$ (better orientation) & 0.597 & 0.615 \\
    \midrule
    Conflict / baseless elements & 5,075 / 9,775 & 1,881 / 3,552 \\
    \bottomrule
  \end{tabular}
\end{adjustbox}
\end{table}

\begin{table}[H]
  \centering
  \caption{Calibration, mean over datasets (lower is better).}
  \label{tab:calib}
  \small
  \begin{adjustbox}{max width=\textwidth}
\begin{tabular}{lcccc}
    \toprule
    & \multicolumn{2}{c}{Links} & \multicolumn{2}{c}{Concepts} \\
    \cmidrule(lr){2-3}\cmidrule(lr){4-5}
    Probability & ECE & Brier & ECE & Brier \\
    \midrule
    \zfmn{} & 0.349 & 0.192 & 0.262 & 0.146 \\
    \zfmn{} (context opinion) & 0.223 & 0.121 & 0.173 & 0.104 \\
    Mean of risks & 0.264 & 0.137 & 0.207 & 0.113 \\
    SelfCheckGPT-NLI & 0.172 & 0.128 & 0.186 & 0.139 \\
    Stacking$^\dagger$ & 0.011 & 0.057 & 0.012 & 0.060 \\
    \bottomrule
  \end{tabular}
\end{adjustbox}
\end{table}

\section{Emotion models}
\label{app:rq4}

\begin{table}[H]
  \centering
  \caption{Emotion models on 71,674 concepts from 1,350 answers (23,273 emotional). Both models also control for task and model; Model 2 uses standard errors clustered by answer.}
  \label{tab:rq4}
  \small
  \begin{adjustbox}{max width=\textwidth}
\begin{tabular}{lrcr}
    \toprule
    Term (per SD) & Coef. & 95\% CI & $p$ \\
    \midrule
    \multicolumn{4}{l}{\emph{Model 1: logit $P$, linear mixed model}} \\
    $|$valence$|$ & $-0.028$ & [$-0.039$, $-0.016$] & $3\times10^{-6}$ \\
    \# emotions & $-0.018$ & [$-0.029$, $-0.006$] & 0.003 \\
    frame valence spread & $+0.513$ & [$+0.503$, $+0.523$] & $<10^{-10}$ \\
    Zipf frequency & $-0.212$ & [$-0.222$, $-0.202$] & $<10^{-10}$ \\
    word length & $-0.150$ & [$-0.160$, $-0.140$] & $<10^{-10}$ \\
    log degree & $+0.211$ & [$+0.201$, $+0.221$] & $<10^{-10}$ \\
    \midrule
    \multicolumn{4}{l}{\emph{Model 2: hallucinated, logistic (odds ratios)}} \\
    risk $1-P$ & 1.826 & [1.725, 1.932] & $<10^{-10}$ \\
    emotional & 0.843 & [0.783, 0.907] & $5\times10^{-6}$ \\
    risk $\times$ emotional & 0.895 & [0.830, 0.964] & 0.003 \\
    \bottomrule
  \end{tabular}
\end{adjustbox}
\end{table}

\section{Released files}
\label{app:release}
\path{ragtruth_fmn_elements.csv.gz} has one row per concept or link with: model, task, RAGTruth response identifier, element type and key (lemma, or two lemmas for a link), number of occurrences, degree, the binary label and RAGTruth hallucination types, EmoAtlas features, and the scores of all methods (\texttt{score\_*}). \path{ragtruth_fmn_answers.csv.gz} has one row per answer with its label and answer-level scores. Scripts rebuild every table and figure from these two files.

\end{document}